\documentclass{amsart}[12pt]
\usepackage{amsmath, amsthm, amscd, amsfonts,}
\usepackage{graphicx,float}

\usepackage{amssymb}

\usepackage{pb-diagram}

\newcommand{\PP}{{\mathbb{P}}}

\newcommand{\url}[1]{{\tt #1}}

\DeclareMathOperator{\dom}{dom}

\def\MQB{{\mathbb{Q}}}

\def\k{\kappa}

\def\a{\alpha}

\newtheorem{theorem}{Theorem}[section]
\newtheorem{lemma}[theorem]{Lemma}

\newtheorem{definition}[theorem]{Definition}

\numberwithin{equation}{section}

\def\MQB{{\mathbb{Q}}}

\def\k{\kappa}

\def\a{\alpha}

\def\rmark{\mbox{$\rm\bf\rule{0.06em}{1.45ex}\kern-0.05em R$}}
\def\pmark{\mbox{$\rm\bf\rule{0.06em}{1.45ex}\kern-0.05em P$}}
\def\nmark{\mbox{$\rm\bf\rule{0.06em}{1.45ex}\kern-0.05em N$}}
\def\vdash{\mbox{$\rm\| \kern-0.13em -$}}
\newcommand{\lusim}[1]{\smash{\underset{\raisebox{1.2pt}[0cm][0cm]{$\sim$}}
{{#1}}}}

\usepackage{pb-diagram}

\def\rmark{\mbox{$\rm\bf\rule{0.06em}{1.45ex}\kern-0.05em R$}}
\def\pmark{\mbox{$\rm\bf\rule{0.06em}{1.45ex}\kern-0.05em P$}}
\def\nmark{\mbox{$\rm\bf\rule{0.06em}{1.45ex}\kern-0.05em N$}}
\def\vdash{\mbox{$\rm\| \kern-0.13em -$}}

\begin{document}

\title[On  a theorem of Woodin]{Power set at $\aleph_\omega$: on  a theorem of Woodin }

\author[M. Golshani.]{Mohammad Golshani}

  \thanks{The author's research was in part supported by a grant from IPM (No. 91030417). He also thanks Radek Honzik for many valuable comments.}
\maketitle

\begin{abstract}
We give Woodin's original proof that if there exists a $(\kappa+2)-$strong cardinal $\kappa,$ then
there is a generic extension of the universe in which $\k=\aleph_\omega,$  $GCH$ holds below $\aleph_\omega$ and $2^{\aleph_\omega}=\aleph_{\omega+2}.$
\end{abstract}

\section{introduction}
One of the central topics in set theory since Cantor has been the study of the power set function $\k \mapsto 2^\k$,  and despite many results which are obtained
about it, determining its behavior is far from being answered completely. In this paper we consider the very special case of determining the power of $2^{\aleph_\omega},$
when $\aleph_\omega$ is a strong limit cardinal, and so just discuss a little about what is known for this case.

The first important results were obtained by Magidor, who proved the consistency of $\aleph_\omega$ is strong limit and $2^{\aleph_\omega}=\aleph_{\omega+k+1},$ where $1<k\leq \omega$, from a supercompact cardinal \cite{magidor1}, and the consistency of $2^{\aleph_\omega}=\aleph_{\omega+2}$ while $GCH$ holding below
$\aleph_\omega$ from a huge cardinal and a supercompact cardinal below it \cite{magidor2}.
In \cite{shelah1}, Shelah improved  Magidor's theorem from \cite{magidor1} by showing the consistency of $\aleph_\omega$ is strong limit and $2^{\aleph_\omega}=\aleph_{\omega+\alpha+1},$ where $\alpha$ is any countable ordinal, still starting from a supercompact cardinal.

In 1980$^{\text{th}}$, Woodin was able to reduce the large cardinal assumptions used by Magidor to the level of strong cardinals, and in particular he proved the following  theorem.
\begin{theorem} (Woodin). Suppose $GCH$ holds and $\kappa$ is a $(\kappa+2)-$strong cardinal. Then there is a generic extension of the universe in which $\k=\aleph_\omega,$  $GCH$ holds below $\aleph_\omega$ and $2^{\aleph_\omega}=\aleph_{\omega+2}.$
\end{theorem}
In \cite{gitik-magidor}, Gitik and Magidor introduced a new method of forcing, called extender based Prikry forcing, and using it, they were able to reduce the large cardinal assumptions used by Magidor and Shelah to the level of strong cardinal.  In particular they proved the following.

\begin{theorem}
Assume $\kappa$ is a $(\kappa+\alpha+1)$-strong cardinal, where $\alpha<\omega_1.$ Then there is a generic extension
$\k=\aleph_\omega,$  $GCH$ holds below $\aleph_\omega$ and $2^{\aleph_\omega}=\aleph_{\omega+\alpha+1}.$
\end{theorem}
In \cite{gitik-merimovich}, Gitik and Merimovich showed that if we allow finite gap at $\aleph_\omega,$ then the continuum function below it can behave arbitrary, in the sense that given any finite $m>1$ and any function  $\phi: \omega \rightarrow \omega$ such that $\phi$ is increasing and $\phi(n)>n$, there is a model of $ZFC$
in which $2^{\aleph_\omega}=\aleph_{\omega+m}$ and for all $n<\omega, 2^{\aleph_n}=\aleph_{\phi(n)}$ (starting from a $(\kappa+m)$-strong cardinal).

Cummings \cite{cummings1} has given a strengthening of Woodin's theorem, by producing a model of $ZFC$ in which $GCH$ holds at all successor cardinals but fails at all limit cardinals, however the proof, which uses Radin forcing is quite complicated, while for Woodin's theorem just a variant of Prikry forcing is sufficient. As there is no explicit proof of Woodin's theorem, in this paper we will sketch a proof of it, which is based on ideas from \cite{cummings1}. We have avoided all
the details, as they all can be found in \cite{cummings1} or \cite{cummings2}.

Let us also mention that by results of Gitik and Mitchel, it is known that some large cardinals at the level of strong cardinals are required for the
results stated above (see \cite{gitik2}, \cite{gitik-mitchell} and \cite{mitchell}).  On the other hand
by a famous result of Shelah, if $\aleph_\omega$ is a strong limit cardinal, then
$2^{\aleph_\omega} < min\{\aleph_{(2^{\aleph_0})^+}, \aleph_{\omega_4}     \}.$

\section{Woodin's result}
In this section we briefly review Woodin's original proof of Theorem 1.1. The proof we present here
is based on ideas from \cite{cummings1}, and we refer to it for details.
Thus assume $GCH$ holds and let $\k$ be a $(\k+2)$-strong cardinal. Let $E$ be a $(\k, \k^{++})$-extender witnessing this, and let
$j: V \rightarrow M \simeq Ult(V, E) \supseteq V_{\k+2}$ be the corresponding elementary embedding with $crit(j)=\k.$
The proof is in several steps.

\subsection{STEP 1.} Factor $j$ through the canonical ultrapower to get the diagram

\begin{center}

\begin{align*}
\begin{diagram}
\node{V}
        \arrow{e,t}{j}
        \arrow{s,l}{i}
        \node{M \simeq Ult(V, E)}
\\
\node{N \simeq Ult(V, E_\k)}
         \arrow{ne,b}{k}
\end{diagram}
\end{align*}
\end{center}
Let
\begin{center}
$\PP^1=\langle \langle \PP^1_\tau: \tau \leq \k+1 \rangle, \langle \lusim{\MQB}^1_\tau: \tau \leq \k \rangle \rangle $
\end{center}
be the  reverse Easton iteration for adding $\tau^{++}$-many Cohen subsets of $\tau^+,$ for each inaccessible $\tau \leq \k.$ So for each $\tau \leq \k, \MQB^1_\tau$ is the trivial forcing, except $\tau \leq \k$ is inaccessible, in which case
 $\Vdash_{\PP^1_{\tau}}$``$ \lusim{\MQB}^1_\tau=\lusim{Add}(\tau^+, \tau^{++}).$

Let
\[
G^1 = \langle  \langle G^1_\tau: \tau \leq \k+1  \rangle, \langle H^1_\tau: \tau \leq \k          \rangle \rangle
\]
be $\PP^1$-generic over $V$ and $V^1=V[G^1].$ Then by standard arguments, there are generic filters $G^1_M, G^1_N \in V^1$ so that the diagram
lifts to the following
\begin{center}

\begin{align*}
\begin{diagram}
\node{V^1}
        \arrow{e,t}{j^1}
        \arrow{s,l}{i^1}
        \node{M^1= M[G^1_M] }
\\
\node{N^1 = N[G^1_N]}
         \arrow{ne,b}{k^1}
\end{diagram}
\end{align*}
\end{center}
The next lemma is essentially proved in \cite{cummings1}.
\begin{lemma}
There exists a  filter $\bar{g} \in V^1$ which is  $i^1(Add(\k, \k^{++})_{V^1} \times Col(\k, \k^+)_{V^1})$-generic over $N^1.$
\end{lemma}

\subsection{STEP 2.} Work in $V^1$. Note that

$\hspace{1.cm}$ $N^1=$the transitive collapse of $\{j^1(f)(\k): f: \k \rightarrow V^1   \},$

$\hspace{1.cm}$ $M^1=$the transitive collapse of $\{j^1(f)(a): f: [\k]^{<\omega} \rightarrow V^1, a \in [\k^{++}]^{<\omega}   \}.$

 Let
 \begin{center}
$\PP^2=\langle \langle \PP^2_\tau: \tau \leq \k \rangle, \langle \lusim{\MQB}^2_\tau: \tau < \k \rangle \rangle $
\end{center}
be the  reverse Easton iteration for adding $\tau^{++}$-many Cohen subsets of $\tau,$ for each inaccessible $\tau \leq \k.$ So for each $\tau < \lambda, \MQB^2_\tau$ is the trivial forcing, except $\tau $ is inaccessible, in which case
 $\Vdash_{\PP^2_{\tau}}$``$ \lusim{\MQB}^2_\tau=\lusim{Add}(\tau, \tau^{++}).$

Let
\[
G^2 = \langle  \langle G^2_\tau: \tau \leq \k  \rangle, \langle H^2_\tau: \tau < \k          \rangle \rangle
\]
be $\PP^2$-generic over $V^1$ and $V^2=V^1[G^2].$
Then for some suitable generic filters in $V^2$, we can lift the diagram one further step and get the following:

\begin{center}
\begin{align*}
\begin{diagram}
\node{V^2}
        \arrow{e,t}{j^2}
        \arrow{s,l}{i^2}
        \node{M^2= M^1[G^2_M] }
\\
\node{N^2 = N^1[G^2_N]}
         \arrow{ne,b}{k^2}
\end{diagram}
\end{align*}
\end{center}
Let $E^2$ denote the $(\k, \k^{++})$-extender derived from $j^2$. We have the following (see \cite{cummings1}).
\begin{lemma}
There exists $F \in V^2$ such that $F$ is $Col(\k^{+3}, < i(\k))_{N^2} \times Col(i(\k), i(\k^+))_{N^2}$-generic over $N^2$. Further $F \in M^2.$
\end{lemma}

\subsection{STEP 3.} In this section we define the main forcing construction. Work in $V^2$. Let $i^2_{\alpha,\beta}: N^2_{\alpha} \rightarrow N^2_\beta$
denote the standard embedding of the $\a^{\text{th}}$ ultrapower into the $\beta^{\text{th}}$ one, and set $\k_\a = i^2_{0,\a}(\k).$
Define
\\
$\PP= [Col(\k^{+3}, < \k_1) \times Col(\k_1, \k^+_1)]_{N^2_{2}},$
\\
$\PP^* = \{f \in N^2_1: \dom(f) \in j^2(E^2_\k): \forall \beta\in dom(f), ~ f(\beta) \in  [Col(\k^{+3}, < \beta) \times Col(\beta, \beta^+)]_{N^2_1}               \}.$

It is not difficult to see that $F$ is $\PP$-generic over $N^2_1$, which gives rise to some $F^*$ which is $\PP^*$-generic over $N^2_1.$
For $\a < \beta$ set
\[
\PP(\alpha, \beta) = Col(\a^{+3}, < \beta) \times Col(\beta, \beta^+).
\]
We now define the main forcing construction.
\begin{definition}
A condition in $\PP^3$ is a finite sequence of the form
\[
\langle \delta_0, P_1, \delta_1, \dots, P_n, \delta_n, H, h               \rangle
\]
where
\begin{enumerate}
\item $\delta_0 < \delta_1 < \dots < \delta_n < \k,$
\item Each $P_k \in \PP(\delta_{k-1}, \delta_k), k \leq n,$
\item $dom(h) \in E^2_\k, dom(h) \subseteq \k \setminus (\delta_n+1),$
\item $h(\beta) \in \PP(\delta_n, \beta),$
\item $dom(H)=[dom(h)]^2,$
\item $H(\a, \beta) \in \PP(\a, \beta),$
\item $i^2_{0,2}(H)(\k, j^2(\k)) \in F.$
\end{enumerate}
\end{definition}
The order relation on $\PP^3$ is defined as follows.
\begin{definition}
Let $p=\langle \delta_0, \dots, P_n, \delta_n, H, h   \rangle$ and $q= \langle \xi_0,  \dots, Q_m, \xi_n, I, i              \rangle$
be two conditions in $\PP^3.$ Then $p \leq q$ if and only if
\begin{enumerate}
\item $n \geq m,$
\item $\forall k \leq m, \delta_k=\xi_k$ and $P_k \leq Q_k,$
\item $\forall k>m,~ \delta_k \in dom(i),$
\item $dom(h) \subseteq dom(i),$
\item  $\forall (\a, \beta) \in dom(H), H(\a, \beta) \leq I(\a, \beta),$
\item $n=m  \Rightarrow h(\beta) \leq i(\beta)$,
\item If $n>m$ then
\begin{itemize}
\item $P_{m+1} \leq i(\delta_{m+1})$,
\item $m+1 < k \leq n \Rightarrow ~ P_k \leq I(\delta_{k-1}, \delta_k),$
\item $h(\beta) \leq I(\delta_n, \beta).$
\end{itemize}
\end{enumerate}
\end{definition}
Let $G^3$ be $\PP^3$-generic over $V^2$
and let $V^3=V^2[G^3].$ Also let
\[
\langle  \delta_n: n<\omega          \rangle
\]
be the Prikry sequence added by $G^3.$
Let us summarize the main properties of the generic extension. The proof is essentially the same as (and in fact simpler than)  the proofs given in \cite{cummings1}, where we also refer to it for more details.
\begin{lemma}
$(1)$ $V^3 \models $``$~\k = \delta_0^{+\omega}$'',

$(2)$ $V^2$ and $V^3$ have the same cardinals $\geq \k,$

$(3)$ In $V^3,~ GCH$ holds in the interval $(\delta_0, \kappa)$ and $2^{\k}=(\k^{++})^{V}.$
\end{lemma}

\subsection{STEP 4.} Force over $V^3$ with $\PP^4=Col(\omega, \delta_0^+),$ and let $G^4$ be $\PP^4$-generic over $V^3$. Also let $V^4=V^3[G^4].$ It is evident that
\[
V^4 \models \text{~``} \aleph_\omega=\k, GCH  ~\text{holds below} ~ \aleph_\omega ~\text{and} ~ 2^{\aleph_\omega} =\aleph_{\omega+2}.
\]

This completes the proof of Theorem 1.1.

School of Mathematics, Institute for Research in Fundamental Sciences (IPM), P.O. Box:
19395-5746, Tehran-Iran.

E-mail address: golshani.m@gmail.com


\begin{thebibliography}{99}

\bibitem{cummings1}  Cummings, James. A model in which GCH holds at successors but fails at limits. Trans. Amer. Math. Soc. 329 (1992), no. 1, 1–-39.

\bibitem{cummings2} Cummings, James. Iterated forcing and elementary embeddings. Handbook of set theory. Vols. 1, 2, 3, 775-–883, Springer, Dordrecht, 2010.




\bibitem{gitik2} Gitik, Moti The strength of the failure of the singular cardinal hypothesis. Ann. Pure Appl. Logic 51 (1991), no. 3, 215-–240.

\bibitem{gitik3} Gitik, Moti, Prikry-type forcings. Handbook of set theory. Vols. 1, 2, 3, 1351–-1447, Springer, Dordrecht, 2010.

\bibitem{gitik-magidor} Gitik, Moti; Magidor, Menachem. The singular cardinal hypothesis revisited. Set theory of the continuum (Berkeley, CA, 1989), 243–-279, Math. Sci. Res. Inst. Publ., 26, Springer, New York, 1992.

\bibitem{gitik-mitchell}   Gitik, Moti; Mitchell, William J. Indiscernible sequences for extenders, and the singular cardinal hypothesis. Ann. Pure Appl. Logic 82 (1996), no. 3, 273–-316.

\bibitem{gitik-merimovich} Gitik, Moti; Merimovich, Carmi, Possible values for $2^{\aleph_n}$ and $2^{\aleph_\omega}$. Ann. Pure Appl. Logic 90 (1997), no. 1--3, 193-–241.

\bibitem{magidor1}  Magidor, Menachem, On the singular cardinals problem. I. Israel J. Math. 28 (1977), no. 1--2, 1-–31.

\bibitem{magidor2}  Magidor, Menachem, On the singular cardinals problem. II. Ann. of Math. (2) 106 (1977), no. 3, 517–-547.


\bibitem{mitchell} Mitchell, William J. The covering lemma. Handbook of set theory. Vols. 1, 2, 3, 1497–-1594, Springer, Dordrecht, 2010.

\bibitem{shelah1} Shelah, Saharon, The singular cardinals problem: independence results. Surveys in set theory, 116,–134, London Math. Soc. Lecture Note Ser., 87, Cambridge Univ. Press, Cambridge, 1983.

\bibitem{shelah2} Shelah, Saharon, Cardinal arithmetic. Oxford Logic Guides, 29. Oxford Science Publications. The Clarendon Press, Oxford University Press, New York, 1994. xxxii+481 pp.





\end{thebibliography}
\end{document}